\documentclass{article}
\usepackage{latexsym}
\usepackage{amssymb}
\usepackage{amsmath}
\newtheorem{theorem}{Theorem}

\newtheorem{lemma}{Lemma}

\newtheorem{example}{Example}
\newtheorem{remark}{Remark}
\def\erre{{\rm I\!R}}
\def\R{{\rm I\!R}}

\def\enne{{\rm I\!N}}
\def\meas{\mathop{\rm meas}}

\begin{document}

\title{On doubly nonlocal fractional elliptic equations}

\author{Giovanni Molica Bisci\footnote{Corresponding author.}\\
Dipartimento P.A.U.\\
          Universit\`a  degli
Studi Mediterranea di Reggio Calabria\\
          Salita Melissari, 89100 Reggio Calabria, Italy\\
          gmolica@unirc.it\\\\
          Du\v{s}an Repov\v{s}\\
          Faculty of Education and Faculty of Mathematics and Physics\\
          University of Ljubljana, POB 2964, Slovenia 1001\\
 dusan.repovs@guest.arnes.si}

\date{}

\maketitle


Scientific Chapters: Mathematics.


\begin{abstract}
This work is devoted to the study of the existence of solutions to nonlocal equations involving the fractional Laplacian.
These equations have a
variational structure and we find a nontrivial solution for them
using the Mountain Pass Theorem.  To make the nonlinear methods work, some
careful analysis of the fractional spaces involved is necessary. In addition, we require rather general assumptions on
the local operator. As far as we know, this result is new and represent a fractional version of a classical theorem obtained working with Laplacian equations.
\end{abstract}

\bigskip

\section{Introduction}
As is well-known, nonlocal Laplacian boundary value problems
model several physical and biological systems where $u$ describes a process which depends on the average of itself, for
example, the population density (see \cite{8,9,10,11,12}). In the vast literature on this subject, we also refer the reader to some interesting results obtained by Autuori and Pucci in \cite{AP1,AP2,AP3}studying Kirchhoff equations by using different approaches.\par
 Very recently, the nonlocal fractional counterpart of Kirchhoff-type problems has been considered (see \cite{caff} and \cite{fv, mio2}). In this order of ideas, we are interested here on the existence of weak solutions for the following (doubly) nonlocal problem:
\begin{equation} \tag{$D_{M,f}$} \label{Nostro}
\left\{
\begin{array}{ll}
M\Big(\displaystyle\int_{\erre^{n}\times\erre^n}\frac{|u(x)-u(y)|^2}{|x-y|^{n+2s}}\,dxdy\Big)(-\Delta)^su
= f(x,u) & \rm in \quad \Omega
\\\\ u=0\,\,{\mbox{ in }} \erre^n\setminus \Omega,\\
\end{array}
\right.
\end{equation}
\noindent where $\Omega$ is a bounded domain in $({\erre}^{n},|\cdot|)$ with smooth
boundary $\partial \Omega$, $s\in (0,1)$ is fixed with $s<n/2$ and $(-\Delta)^s$ is the fractional Laplace operator,
which (up to normalization factors) can be defined as
$$
(-\Delta)^s u(x):=
-
\int_{\erre^n}\frac{u(x+y)+u(x-y)-2u(x)}{|y|^{n+2s}}\,dy,
\,\,\,\,\, x\in \erre^n.
$$
Further, $f:\bar\Omega\times \erre\rightarrow \erre$ and $M:[0,+\infty)\to [0,+\infty)$ are suitable continuous maps.\par
In our context, problem \eqref{Nostro} is studied by exploiting classical variational methods. More precisely, we apply the celebrated Mountain Pass Theorem (abbreviated MPT) to this kind of equations motivated by the current literature where the MPT has been intensively applied to find solutions to
quasilinear elliptic equations (see \cite{ar,puccirad, rabinowitz,struwe}).\par
\indent Technically, this approach is realizable by checking that the associated energy functional satisfies the usual compactness Palais-Smale condition in a suitable variational setting developed by Servadei and Valdinoci (see \cite{svmountain}). Indeed, the nonlocal analysis which we perform in this paper
in order to use the Mountain Pass Theorem is quite general and it was successfully exploited for other goals
in several recent contributions (see \cite{mio1,Molica, Molica1, Molica2, Molica3, sY, svmountain,svlinking,servadeivaldinociBN} and \cite{valpal} for various properties on the fractional Sobolev space setting).\par
This functional analytical context was inspired by (but is not equivalent
to) the fractional Sobolev spaces, in order to correctly encode
the Dirichlet boundary datum in the variational formulation.\par
In our context, to avoid some additional technical difficulties due to the presence of the term
$$M\Bigg(\displaystyle\int_{\erre^{n}\times\erre^n}\frac{|u(x)-u(y)|^2}{|x-y|^{n+2s}}\,dx\,dy\Bigg),$$
we impose some restrictions on the behavior of the continuous map $M$.\par
\indent
 More precisely, we require that
there exists a constant $m_0$ such that:
\begin{itemize}
\item[$(C_M^1)$] $0< m_0\leq M(t),\,\,\,\forall\, t\in [0,+\infty)$.
\end{itemize}

\noindent In addition to the above hypothesis, we assume that:
 \begin{itemize}
 \item[$(C_M^2)$] \textit{There exists} $t_0\geq 0$ \textit{such that}
 $$
\displaystyle{\displaystyle \widehat{M}(t)}\geq t M(t),
$$
\textit{for every} $t\in [t_0,+\infty)$, \textit{where} $\displaystyle
\widehat{M}(t):=\int_0^t M(s)ds.
$
\end{itemize}
\indent The above conditions ensure, as proved in Lemma \ref{lemmaSer}, that the potential $\widehat{M}$
has a sublinear growth. Under the previous assumptions, by imposing conditions on the nonlinear part $f$ (among others, the Ambrosetti-Rabinowitz relation) we prove the existence of at
least one nontrivial weak solution to problem \eqref{Nostro}, see Theorem \ref{Esistenza}.\par
 This result is related to \cite[Theorem 2]{svmountain} where the authors studied a local problem involving a general
integro-differential operator of fractional type (see Remark \ref{GO}) whose prototype is
\begin{equation} \tag{$D_{f}$} \label{SV}
\left\{
\begin{array}{ll}
(-\Delta)^su
= f(x,u) & \rm in \quad \Omega
\\u=0 & \mbox{in}\,\,\,\erre^n\setminus\Omega,
\end{array}
\right.
\end{equation}

\indent We just observe that in our context, in contrast with the cited result, we don't require that
$$
\lim_{t\rightarrow 0}\frac{f(x,t)}{t}=0,
$$
uniformly with respect to $x\in\bar\Omega$ (see condition \eqref{atzero} in Theorem \ref{Esistenza} and Remark \ref{finale} below).
\par

Moreover, our main theorem extends to the nonlocal setting a result, already known in the literature for Kirchhoff-type problems obtained by Alves, Corr\^{e}a and Ma \cite[Theorem 3]{ACM}. We just point out that $M$, in the original meaning for Kirchhoff equation, is an increasing function, thus condition $(C_M^2)$ is clearly violated.\par
 We mention, for completeness, that the existence and multiplicity of solutions for elliptic equations in $\erre^n$, driven by a nonlocal integro-differential operator, whose standard prototype is the fractional Laplacian, have been studied by Autuori and Pucci in \cite{AP} (this work is related to the results on general quasilinear elliptic problems given in \cite{AP0}). See also the relevant contributions \cite{p3,p1,p2} where Kirchhoff-type problems have been studied by using different methods and approaches.\par
\indent The plan of the paper is as follows. Section 2 is devoted to our abstract framework and preliminaries. Successively, in Sections 3 we give the main result (see Theorem \ref{Esistenza}). Finally, a concrete example of an application is presented in the last part of the paper (see Example \ref{esempio0}).
\section{Abstract Framework}\label{section2}
 This section is devoted to the notations used throughout the paper. We
also list some preliminary results which will be useful in the sequel.\par
Let $H^s(\erre^n)$ be the usual fractional Sobolev space endowed with the norm (the
so-called \emph{Gagliardo norm})
$$
\|g\|_{H^s(\erre^n)}=\|g\|_{L^2(\erre^n)}+
\Big(\int_{\erre^n\times\erre^n}\frac{\,\,\,|g(x)-g(y)|^2}{|x-y|^{n+2s}}\,dxdy\Big)^{1/2}\,.
$$
and defined as the linear space of functions $g \in L^2(\erre^n)$ such that
$$\mbox{the map}\,\,\,
(x,y)\mapsto \frac{g(x)-g(y)}{|x-y|^{n/2+s}}\,\,\, \mbox{is in}\,\,\,
L^2(\erre^n\times \erre^n, dxdy\big).$$
\indent Let us consider the subspace~$X_0\subset H^s(\erre^n)$ given by
$$
X_0:=\big\{g\in H^s(\erre^n) : g=0\,\, \mbox{a.e. in}\,\,
\erre^n\setminus
\Omega\big\}.
$$
Of course, the space~$X_0$ is non-empty, since $C^2_0 (\Omega)\subseteq X_0$ by \cite[Lemma~11]{sv} and it depends on the set $\Omega$.\par
Moreover, by \cite[Lemma~6]{svmountain} and the fact that any function $v\in X_0$ is such that $v=0$ a.e. in $\erre^n\setminus \Omega$, we can take in the sequel
$$
X_0\ni v\mapsto \|v\|_{X_0}:=\left(\int_{Q}\frac{|v(x)-v(y)|^2}{|x-y|^{n+2s}}\,dxdy\right)^{1/2}
$$
as norm on $X_0$, where $Q:=(\erre^n\times\erre^n)\setminus ({\mathcal{C}}\Omega\times{\mathcal{C}}\Omega)$, and ${\mathcal{C}}\Omega:=\erre^n\setminus\Omega$.\par
Also $\left(X_0, \|\cdot\|_{X_0}\right)$ is a Hilbert space (for this see \cite[Lemma~7]{svmountain}), with the scalar product
$$
\langle u,v\rangle_{X_0}:=\int_{Q}
\frac{\big( u(x)-u(y)\big) \big( v(x)-v(y)\big)}{|x-y|^{n+2s}}\,dxdy.
$$
\indent We recall that in \cite[Lemma~8]{svmountain} and \cite[Lemma~9]{servadeivaldinociBN} the authors proved that the embedding $j:X_0\hookrightarrow
L^{\nu}(\erre^n)$ is continuous for any $\nu\in [1,2^*]$, while it is
compact whenever $\nu\in [1,2^*)$, where $2^*:={2n}/{(n-2s)}$ denotes the Sobolev fractional exponent.\par
 Hence, for any $\nu\in [1, 2^*)$, there exists $c_\nu>0$ such that
$$
\|v\|_{L^{\nu}(\erre^n)}\leq c_{\nu} \|v\|_{X_0},
$$
for every $v\in X_0$.\par
\indent In the sequel, we will denote by $\lambda_{1,s}$ the first (simple and positive) eigenvalue of
the operator~$(-\Delta)^s$ with homogeneous Dirichlet boundary
data, namely the first eigenvalue of the problem
$$\left\{\begin{array}{ll}
(-\Delta)^s u=\lambda u & \mbox{in}\,\,\,\Omega\\
u=0 & \mbox{in}\,\,\,\erre^n\setminus\Omega,
\end{array}\right.$$
that can be characterized as follows
$$
\lambda_{1,s}=\min\left\{\frac{\displaystyle\int_{Q}\frac{|u(x)-u(y)|^2}{|x-y|^{n+2s}}\,dx\,dy}{\displaystyle\int_{\Omega}|
u(x)|^{2}dx}:u\in X_0\setminus\{0_{X_0}\}\right\}.
$$

\indent For the existence and the basic properties of this eigenvalue we
refer to \cite[Proposition~9 and Appendix~A]{svlinking}, where a
spectral theory for general integrodifferential nonlocal operators
was developed. Further properties can be also found in \cite{sY}.

Finally, for the sake of completeness, we recall that a $C^1$-functional $J:E\to\R$, where $E$ is a real Banach
space with topological dual $E^*$, satisfies the \textit{Palais-Smale condition at level} $\mu\in\R$, (abbreviated $\textrm{(PS)}_{\mu}$) when:
\begin{itemize}
\item[$\textrm{(PS)}_{\mu}$] {\it Every sequence $\{u_j\}$ in $E$ such that
$$
J(u_j)\to \mu, \quad{\rm and}\quad \|J'(u_j)\|_{E^*}\to0,
$$
as $j\rightarrow \infty$, possesses a convergent subsequence.}
\end{itemize}
We say that $J$ satisfies the \textit{Palais-Smale condition} (abbreviated $\textrm{(PS)}$) if $\textrm{(PS)}_{\mu}$ holds for every $\mu\in \erre$.\par
With the above notation, our main tool is the classical MPT:

\begin{theorem}\label{MPT}
Let $(E,\|\cdot\|_E)$ be a real Banach space and let $J\in C^1(E;\erre)$ be such that $J(0_E)=0$ and it satisfies the $(\rm PS)$ condition. Suppose that$:$
 \begin{itemize}
 \item[$(I_1)$] There exist constants $\rho,\alpha>0$ such that $J(u)\geq \alpha$ if $\|u\|_E=\rho$.
  \item[$(I_2)$] There exists $e\in E$ with $\|e\|_E> \rho$ such that $J(e)\leq 0$.
 \end{itemize}
 Then $J$ possesses a critical value $c\geq \alpha$, which can be characterized as
 $$
 c:=\inf_{\gamma\in \Gamma}\max_{u\in \gamma([0,1])}J(u),
 $$
 where
 $$
 \Gamma:=\{\gamma\in C([0,1];E): \gamma(0)=0\,\,\wedge\,\, \gamma(1)=e\}.
 $$
\end{theorem}

See \cite[p. 7; Theorem 2.2]{rabinowitz}.\par
\smallskip
 \noindent We cite the monograph \cite{k2} as general reference for the variational setting adopted in this paper.
\section{The main result}
Our main result is as follows.
\begin{theorem}\label{Esistenza}
Let us assume that $M:[0,+\infty)\to [0,+\infty)$ is a continuous map such that conditions $(C_M^1)$ and $(C_M^2)$ hold. Further, require that $f:\bar\Omega\times\erre\rightarrow \erre$ is a continuous function which satisfies the following requirements$:$
\begin{itemize}
\item[$\rm h_1)$] The subcritical growth condition$:$
$$
 |f(x,t)|\leq c(1+|t|^{q-1}),\quad(\forall\, x\in \bar\Omega,\,\forall\,t\in \erre)
$$
where $c>0$ and $2<q<2^{*}$$;$
\item[$\rm h_2)$] The Ambrosetti-Rabinowitz $($abbreviated $(\rm AR)$$)$ condition$:$ $$F(x,\xi):=\displaystyle\int_0^\xi f(x,t)dt$$ is $\theta$-superhomogeneous at infinity, i.e. there exists $t^{\star}>0$ such that
$$
0<\theta F(x,\xi)\leq f(x,\xi)\xi,\quad(\forall\,x\in \bar\Omega,\,\forall\, |\xi|\geq t_{\star})
$$
    where $\theta>2$.
\end{itemize}
We also assume that
\begin{equation}\label{atzero}
 \displaystyle\limsup_{t\rightarrow 0}\frac{f(x,t)}{t}\leq \lambda,
 \end{equation}
 uniformly for $x\in\bar \Omega$, where
 $$
 \lambda< m_0\lambda_{1,s}.
 $$
Then the nonlocal problem
\begin{equation} \tag{$D_{M,f}$} \label{Nostro2}
\left\{
\begin{array}{ll}
M\Big(\displaystyle\int_{Q}\frac{|u(x)-u(y)|^2}{|x-y|^{n+2s}}\,dxdy\Big)(-\Delta)^su
= f(x,u) & \rm in \quad \Omega
\\\\ u=0\,\,{\mbox{\rm in }} \erre^n\setminus \Omega,\\
\end{array}
\right.
\end{equation}
\noindent has at least one nontrivial weak solution.
\end{theorem}

We recall that a \textit{weak solution} of problem \eqref{Nostro} is a function $u\in X_0$ such that
$$
\begin{array}{l} {\displaystyle M\Big(\displaystyle\int_{Q}\frac{|u(x)-u(y)|^2}{|x-y|^{n+2s}}dxdy\Big)\int_{Q}
\frac{(u(x)-u(y))(\varphi(x)-\varphi(y))}{|x-y|^{n+2s}} dxdy}\\
\\
\displaystyle \qquad\qquad\qquad\qquad\qquad\quad\quad= {\int_\Omega f(x,u(x))\varphi(x)dx},\,\,\,\,\,\,\, \forall\,\, \varphi \in X_0.
\end{array}
$$
\subsection{Some remarks on our assumptions}
The validity of the next lemma will be crucial in the sequel.
 \begin{lemma}\label{lemmaSer}
 Suppose that conditions $(C_M^1)$ and $(C_M^2)$ hold. Then, there are two positive constants $m_1$ and $m_2$ such that
 \begin{equation}\label{cresce}
 \widehat{M}(t)\leq m_1t +m_2,
 \end{equation}
 for every $t\in [0,+\infty)$.
\end{lemma}
 \textit{Proof.}
 Let $t_1>t_0$, where $t_0$ appears in hypothesis $(C_M^2)$. By our assumptions we easily have
 $$
 \frac{M(t)}{\widehat{M}(t)}\leq \frac{1}{t},
 $$
 for every $t\in [t_1,+\infty)$. Integrating the above relation, we obtain
 $$
\int_{t_1}^{t}\frac{M(s)}{\widehat{M}(s)}ds= \log \frac{\widehat{M}(t)}{{\widehat{M}}(t_1)}\leq \log \frac{t}{t_1},
 $$
 for every $t\in]t_1,+\infty)$. Thus
 $$
 \widehat{M}(t)\leq \frac{\widehat{M}(t_1)}{t_1}t,
 $$
 for every $t\in ]t_1,+\infty)$. Hence the growth condition \eqref{cresce} holds by taking, for instance, $m_1:=\displaystyle\frac{\widehat{M}(t_1)}{t_1}$ and $m_2:=\displaystyle\max_{t\in [0,t_1]}\widehat{M}(t)$. The proof is complete.
 
 \smallskip
 
 Owing to conditions $(C_M^1)$ and $(C_M^2)$, by Lemma \ref{lemmaSer} one gets the following inequalities:
 \begin{itemize}
 \item[$(\widehat{C}_M)$]
 $\displaystyle m_0\frac{\|u\|^{2}_{X_0}}{2}\leq \frac{1}{2}\widehat{M}\left(\displaystyle\int_{Q}\frac{|u(x)-u(y)|^2}{|x-y|^{n+2s}}\,dxdy\right)\leq {\displaystyle m_1\frac{\|u\|^{2}_{X_0}}{2} +\frac{m_2}{2}}$
 \end{itemize}
 for every $u\in X_0$.
  \begin{remark}\rm{
 We observe that $\rm h_2)$ implies that
 $$
 F(x,\tau \xi)\geq F(x,\xi)\tau^\theta,
 $$
 for every $x\in \bar\Omega$, $|\xi|\geq t_{\star}$ and $\tau\geq 1$. Indeed, for $\tau=1$, clearly the equality holds. Otherwise, fix $|\xi|\geq t_{\star}$ and define $g(x,\tau):=F(x,\tau\xi)$, for every $x\in\Bar\Omega$ and $\tau\in ]1,+\infty)$. By $(\rm AR)$ condition it follows that
 $$
 \frac{g'(x,\tau)}{g(x,\tau)}\geq \frac{\theta}{\tau},
 $$
 for every $x\in \bar\Omega$ and $\tau> 1$. By integrating in $]1,\tau]$ we get that
 $$
\int_{1}^{\tau}\frac{g'(x,s)}{g(x,s)}ds= \log \frac{g(x,\tau)}{g(x,1)}\geq \log \tau^{\theta}.
 $$
In conclusion, since for every $x\in \bar\Omega$, $|\xi|\geq t_{\star}$ and $\tau> 1$ one has
 \begin{eqnarray*}
F(x,\tau\xi) &=:& g(x,\tau)\\\nonumber
&\geq& g(x,1)\tau^{\theta}\\\nonumber
&=& F(x,\xi)\tau^\theta,\nonumber
\end{eqnarray*}
 the claim is verified.}
 \end{remark}
\subsection{Proof of Theorem \ref{Esistenza}} Set
$$
\Phi(u):=\frac12\widehat{M}\left(\displaystyle\int_{Q}\frac{|u(x)-u(y)|^2}{|x-y|^{n+2s}}\,dxdy\right)
$$
for every $u\in X_0$.\par
Under the assumptions of Theorem \ref{Esistenza}, we define the $C^1$-functional
 $$
 J(u):=\Phi(u)-\int_\Omega F(x,u(x))dx,\,\,\,\ (\forall\,u\in X_0)
 $$
 where
 $$
 F(x,u(x)):=\int_0^{u(x)}f(x,s)ds,
 $$
whose critical points are the weak solutions of problem \eqref{Nostro}.\par
 In order to prove our result, we apply Theorem \ref{MPT} to this functional. In the next three lemmas we shall verify the Mountain Pass Theorem conditions.
\begin{lemma}\label{lemma1}
 Every Palais-Smale sequence for the functional $J$
 is bounded in $X_0$.
\end{lemma}
\textit{Proof.}
Let $\{u_j\}\subset X_0$ be a Palais-Smale sequence i.e.
\begin{equation}\label{E1}
J(u_j)\rightarrow \mu,
\end{equation}
for $\mu\in \erre$ and
\begin{equation}\label{E2}
\quad \|J'(u_j)\|_{X_0^*}=\sup\Big\{ \big|\langle\,J'(u_j),\varphi\,\rangle \big|\,: \;
\varphi\in
X_0\,, \|\varphi\|_{X_0}=1\Big\}\to0,
\end{equation}
as $j\rightarrow \infty$.\par
\indent We argue by contradiction. So, suppose to the contrary that the conclusion were not true. Passing to a subsequence if necessary, we may assume that
$$
\|u_j\|_{X_0}\rightarrow +\infty,
$$
as $j\rightarrow \infty$.\par
\indent  By conditions $(C_M^1)$ and $(C_M^2)$, it follows that there exists $j_0\in \enne$ such that
  \begin{eqnarray*}
J(u_j) &-& \frac{\langle J'(u_j),u_j\rangle }{\theta}\\\nonumber
&=& M(\|u_j\|^2_{X_0})\left[\frac{\widehat{M}(\|u_j\|^2_{X_0})}{2 M(\|u_j\|^2_{X_0})}-\frac{\|u_j\|^2_{X_0}}{\theta}\right]\\\nonumber
               &+& \int_{\Omega}\left[\frac{f(x,u_j(x))u_j(x)}{\theta}-F(x,u_j(x))\right]dx,\\\nonumber
               &\geq& m_0\left(\frac{1}{2}-\frac{1}{\theta}\right)\|u_j\|^2_{X_0}\\\nonumber
               &+& \int_{\Omega}\left[\frac{f(x,u_j(x))u_j(x)}{\theta}-F(x,u_j(x))\right]dx,\nonumber
\end{eqnarray*}
for every $j\geq j_0$.\par
 Thus
\begin{eqnarray*}
m_0\left(\frac{\theta-2}{2\theta}\right)\|u_j\|_{X_0}^2&\leq& J(u_j)-\frac{\langle J'(u_j),u_j\rangle }{\theta}\\\nonumber
               &-& \int_{|u_j(x)|>t_{\star}}\left[\frac{f(x,u_j(x))u_j(x)}{\theta}-F(x,u_j(x))\right]dx,\\\nonumber
               &+& M \meas(\Omega),\quad\, \forall\, j\geq j_0,
\end{eqnarray*}
where $``\meas(\Omega)"$ denotes the standard Lebesgue measure of $\Omega$ and
$$
M:=\sup\left\{\left|\frac{f(x,t)t}{\theta}-F(x,t)\right|:x\in\bar\Omega, |t|\leq t_{\star}\right\}.
$$

\indent Now, we observe that, the $(\rm AR)$ condition yields
$$
\int_{|u_j(x)|>t_{\star}}\left[\frac{f(x,u_j(x))u_j(x)}{\theta}-F(x,u_j(x))\right]dx\geq 0.
$$
So, we deduce that
\begin{eqnarray*}
m_0\left(\frac{\theta-2}{2\theta}\right)\|u_j\|_{X_0}^2\leq J(u_j)-\frac{\langle J'(u_j),u_j\rangle }{\theta}+M \meas(\Omega),
\end{eqnarray*}
for every $j\geq j_0$.\par
Then, for every $j\geq j_0$ one has
\begin{eqnarray*}
C\|u_j\|_{X_0}^2\leq J(u_j)+{\theta}{\|J'(u_j)\|_{X_0^*}\|u_j\|_{X_0}}+M \meas(\Omega),
\end{eqnarray*}
where $C:=\displaystyle m_0\left(\frac{\theta-2}{2\theta}\right)>0$.\par
\noindent In conclusion, dividing by $\|u_j\|_{X_0}$ and letting $j\rightarrow \infty$, we obtain a contradiction.

\smallskip

The above lemma implies that the $C^1$-functional
$J$ satisfies the Palais-Smale condition as proved in the next result.
\begin{lemma}\label{lemma2}
The functional $J$
satisfies the compactness $(\rm PS)$ condition.
\end{lemma}

\noindent \textit{Proof.} Let $\{u_{j}\}\subset X_0$ be a Palais-Smale
sequence. By Lemma \ref{lemma1}, the sequence $\{u_j\}$ is necessarily bounded in $X_0$. Since $X_0$ is reflexive, we can extract a subsequence
which for simplicity we shall call again $\{u_{j}\}$, such that
$u_{j}\rightharpoonup u_\infty$ in $X_0$. This means that
\begin{equation}\label{convergenze0}
\begin{aligned}
 & \int_Q\frac{\big(u_j(x)-u_j(y)\big)\big(\varphi(x)-\varphi(y)\big)}{|x-y|^{n+2s}}dxdy \to \\
& \qquad \qquad \qquad
\int_Q\frac{\big(u_\infty(x)-u_\infty(y)\big)\big(\varphi(x)-\varphi(y)\big)}{|x-y|^{n+2s}}dxdy,
\end{aligned}
\end{equation}
for any $\varphi\in X_0$,
as $j\to +\infty$.\par
 We will prove that ${u_{j}}$ strongly converges to $u_\infty\in X_0$.
Exploiting the derivative $J'(u_j)(u_j-u_\infty)$, we obtain
\begin{eqnarray}\label{jj}
\langle a(u_{j}),u_{j}-u_\infty\rangle &=& \langle J'(u_j),u_j-u_\infty\rangle\\
               &+ & \int_{\Omega}f(x,u_j(x))(u_j-u_\infty)(x)dx,\nonumber
\end{eqnarray}
where we set
\begin{eqnarray*}
\langle a(u_{j}),u_{j}-u_\infty\rangle &:=& \Bigg(\int_{Q} \frac{\left|u_j(x)-u_j(y)\right|^{2}}{|x-y|^{n+2s}}dxdy\\\nonumber
               &-& \int_{Q} \frac{\big(u_j(x)-u_j(y)\big)\big(u_\infty(x)-u_\infty(y)\big)}{|x-y|^{n+2s}}dxdy\Bigg)\\\nonumber
               &\times& M\Bigg(\displaystyle\int_{Q}\frac{|u_j(x)-u_j(y)|^2}{|x-y|^{n+2s}}dxdy\Bigg).
\end{eqnarray*}
\indent Since $\|J'(u_j)\|_{X_0^{*}}\to0$ and the sequence $\{u_j-u_\infty\}$ is bounded in
$X_0$, taking into account that $|\langle
J'(u_j),u_j-u_\infty\rangle|\leq\|J'(u_j)\|_{X_0^{*}}\|u_j-u_\infty\|_{X_0}$, one gets
\begin{eqnarray}\label{j2}
\langle J'(u_j),u_j-u_\infty\rangle\to0,
\end{eqnarray}
as $j\rightarrow \infty$.\par
 \indent Now observe that, by condition $\textrm{h}_1)$, it follows that
\begin{align*}
&\int_{\Omega}|f(x,u_j(x))||u_j(x)-u_\infty(x)|dx \\
&\leq c\left(\int_{\Omega}|u_j(x)-u_\infty(x)|dx
+ \int_{\Omega}|u_j(x)|^{q-1}|u_j(x)-u_\infty(x)|dx\right) \\
&\leq c((\meas(\Omega))^{1/q'}+\|u_j\|_{L^{q}(\Omega)}^{q-1}) \|u_j-u_\infty\|_{L^q(\Omega)},
\end{align*}
where, as usual, $q'$ denotes the conjugate of $q$.\par
Since the embedding
$X_0\hookrightarrow L^q(\Omega)$ is compact, clearly $u_j\to u_\infty$ strongly in $L^q(\Omega)$.
So we obtain
\begin{eqnarray}\label{j3}
\int_{\Omega}|f(x,u_j(x))||u_j(x)-u_\infty(x)|dx\to0.
\end{eqnarray}
 \indent By \eqref{jj}, relations \eqref{j2} and \eqref{j3} yield
\begin{eqnarray}\label{fin}
\langle a(u_{j}),u_{j}-u_\infty\rangle
\rightarrow 0,
\end{eqnarray}
when $j\rightarrow \infty$.\par
Bearing in mind condition $(C_M^1)$ one obtains
\begin{eqnarray}\label{finemu}
0<m_0\leq M\Big(\displaystyle\int_{Q}\frac{|u_j(x)-u_j(y)|^2}{|x-y|^{n+2s}}\,dxdy\Big),
\end{eqnarray}
for every $j\in\enne$.\par
Hence by \eqref{finemu} and \eqref{fin} we can write
\begin{eqnarray}\label{fin2}
&& \int_{Q} \frac{\left|u_j(x)-u_j(y)\right|^{2}}{|x-y|^{n+2s}}dxdy-\\\nonumber
               && \quad\quad\int_{Q} \frac{\big(u_j(x)-u_j(y)\big)\big(u_\infty(x)-u_\infty(y)\big)}{|x-y|^{n+2s}}dxdy\rightarrow 0,\nonumber
\end{eqnarray}
when $j\rightarrow \infty$.\par
Thus, by \eqref{fin2} and \eqref{convergenze0} it follows that
\begin{eqnarray*}
&& \limsup_{j\rightarrow \infty}\int_{Q} \frac{\left|u_j(x)-u_j(y)\right|^{2}}{|x-y|^{n+2s}}dxdy\\\nonumber
               && \quad\quad\quad\quad\quad\quad=\int_{Q} \frac{\left|u_\infty(x)-u_\infty(y)\right|^{2}}{|x-y|^{n+2s}}dxdy.\nonumber
\end{eqnarray*}
\indent In conclusion, thanks to \cite[Proposition III.30]{brezis}, $u_j\rightarrow u_\infty$ in $X_0$. The proof is thus complete.

\begin{lemma}\label{lemma3}
The functional $J$ has the geometry of the Mountain Pass Theorem. More precisely$:$
\begin{itemize}
\item[1.] There exists $r>0$ such that
$$
\inf_{\|u\|_{X_0}=r}J(u)>0.
$$
\item[2.] For some $u_0\in X_0$ one has
$$
J(\tau u_0)\rightarrow -\infty,
$$
as $\tau\rightarrow +\infty$.
\end{itemize}
\end{lemma}
\textit{Proof.}
1. We choose $\varepsilon>0$ small enough, so that it satisfies
$$
m_0>\frac{\lambda+\varepsilon}{\lambda_{1,s}}.
$$

\noindent By condition \eqref{atzero} there exists $\delta_\varepsilon>0$ such that
$$
\frac{f(x,t)}{t}\leq \lambda+\varepsilon,
$$
for every $x\in \bar\Omega$ and $|t|\leq\delta_\varepsilon$.\par
 Hence, one has
$$
F(x,\xi)\leq \frac{\lambda+\varepsilon}{2}|\xi|^2,
$$
for every $|\xi|\leq \delta_\varepsilon$.\par
\indent As a consequence of the above inequality, using hypotheses $\rm h_1)$, the Sobolev embedding $X_0\hookrightarrow L^q(\Omega)$ and $(\widehat{C}_M)$, we can write:
\begin{eqnarray*}
J(u) &\geq& \displaystyle\frac{m_0}{2}\|u\|^{2}_{X_0}-\int_{|u(x)|\leq \delta_\varepsilon}\frac{\lambda+\varepsilon}{2}|u(x)|^2dx\\\nonumber
&-& C\int_{|u(x)|> \delta_\varepsilon}|u(x)|^qdx\\\nonumber
               &\geq& \displaystyle\frac{m_0}{2}\|u\|^{2}_{X_0}-\frac{\lambda+\varepsilon}{2\lambda_{1,s}}\|u\|^{2}_{X_0}-C\|u\|^{q}_{X_0},\nonumber
\end{eqnarray*}
for a suitable positive constant $C$.\par
Now, set $r:=\|u\|^2_{X_0}$ and observe that for $r>0$ small enough, we have
$$
\frac{1}{2}\left(m_0-\frac{\lambda+\varepsilon}{\lambda_{1,s}}\right)r-Cr^{q/2}>0,
$$
bearing in mind that $q>2$. Hence
$$
\inf_{\|u\|_{X_0}=r}J(u)>0.
$$
\noindent 2. Let us choose an element $u_0\in X_0$ such that
$$
\meas\left(\{x\in \Omega:u_0(x)\geq t_{\star}\}\right)>0.
$$

\indent $F(x,\xi)$ being a $\theta$-superhomogeneous function if $|\xi|\geq t_{\star}$, for $\tau>1$, we have that
\begin{eqnarray*}
J(\tau u_0) &\leq& \displaystyle{\frac{m_1}{2}\|\tau u_0\|^{2}_{X_0} +\frac{m_2}{2}}-\int_\Omega F(x,\tau u_0(x))dx\\\nonumber
&\leq& \displaystyle m_1\frac{\| u_0\|^{2}_{X_0}}{2}\tau^2 -\tau^\theta\int_{|u_0(x)|\geq t_\star}F(x,u_0(x))dx\\\nonumber
               && +\frac{m_2}{2}+M\meas(\Omega),\nonumber
\end{eqnarray*}
where
$$
M:=\sup\left\{\left|F(x,\xi)\right|:x\in\bar\Omega, |\xi|\leq t_{\star}\right\}.
$$

\indent Thus the $(\rm AR)$ condition implies that
$$
J(\tau u_0)\rightarrow -\infty,
$$
as $\tau\rightarrow +\infty$. This completes the proof.
\section{An example of an application}
 In this section we present a simple application of our main result.
\begin{example}\label{esempio0}\rm
 Consider the following nonlocal problem:
\begin{equation} \tag{$D_{M}$} \label{Nostrom}
\left\{
\begin{array}{ll}
M\Big(\displaystyle\int_{Q}\frac{|u(x)-u(y)|^2}{|x-y|^{n+2s}}\,dxdy\Big)(-\Delta)^su
= u^3+u^4 & \rm in \quad \Omega
\\\\ u=0\,\,{\mbox{ in }} \erre^n\setminus \Omega,\\
\end{array}
\right.
\end{equation}
where
$$
\small M\left(\displaystyle\int_{Q}\frac{|u(x)-u(y)|^2}{|x-y|^{n+2s}}\,dxdy\right):=2+\frac{\sin \left(\displaystyle\int_{Q}\frac{|u(x)-u(y)|^2}{|x-y|^{n+2s}}\,dxdy\right)}{1+\left(\displaystyle\int_{Q}\frac{|u(x)-u(y)|^2}{|x-y|^{n+2s}}\,dxdy\right)^2}.
$$

\indent By virtue of Theorem \ref{Esistenza}, problem \eqref{Nostrom} admits one nontrivial weak solution. Indeed, a direct computation ensures that the continuous function
$$\displaystyle M(t):=2+\frac{\sin t}{1+t^2},\quad (\forall\,t\geq 0)$$ satisfies conditions $(C_M^1)$ and $(C_M^2)$.\par
 Further, the function $f:\erre \rightarrow \erre$ given by $$f(t):=t^{3}+t^4,\quad (\forall\,t\in \erre)$$ satisfies all the hypotheses of Theorems \ref{Esistenza}.
\end{example}

\begin{remark}\rm{
We note that in Example \ref{esempio0} condition $(C_M^2)$ is satisfied for every $t\geq t_0$, where $t_0$ is the unique positive solution of the following real equation
$$
\int_0^t\left(2+\frac{\sin s}{1+s^2}\right)ds-t\left(2+\frac{\sin t}{1+t^2}\right)=0.
$$}
\end{remark}
\begin{remark}\label{GO}\rm{
We just observe that Theorem \ref{Esistenza} can be proved for a more
general class of nonlocal problems of the form
\begin{equation*}
\left\{
\begin{array}{ll}
-M\Big(\displaystyle\int_{\erre^n\times\erre^n}|v(x)-v(y)|^2K(x-y)dxdy\Big)\mathcal L_Ku
= f(x,u) & \rm in \quad \Omega
\\u=0\,\,{\mbox{ in }} \erre^n\setminus \Omega,\\
\end{array}
\right.
\end{equation*}
where $\mathcal L_K$ is defined as follows:
\begin{equation*}
\mathcal L_Ku(x):=
\int_{\erre^n}\Big(u(x+y)+u(x-y)-2u(x)\Big)K(y)\,dy,
\,\,\, (x\in \erre^n)
\end{equation*}
and $K:\erre^n\setminus\{0\}\rightarrow(0,+\infty)$ is a function with the properties that:
\begin{itemize}
\item[1.] $\gamma K\in L^1(\erre^n)$, \textit{where} $\gamma(x):=\min \{|x|^2, 1\}$;
\item[2.] \textit{There exists} $\beta>0$ \textit{and} $s\in (0,1)$
\textit{such that} $$K(x)\geq \beta |x|^{-(n+2s)},$$
{\textit{for any}} $x\in \erre^n \setminus\{0\}$.
\end{itemize}
In this case we look for (weak) solutions $u\in X_0$, where
$$X_0:=\big\{g\in X : g=0\,\, \mbox{a.e. in}\,\,
\erre^n\setminus
\Omega\big\}.$$
Here $X$ denotes the linear space of Lebesgue
measurable functions from $\erre^n$ to $\erre$ such that the restriction
to $\Omega$ of any function $g$ in $X$ belongs to $L^2(\Omega)$ and
$$
((x,y)\mapsto (g(x)-g(y))\sqrt{K(x-y)})\in
L^2\big((\erre^n\times\erre^n) \setminus ({\mathcal C}\Omega\times
{\mathcal C}\Omega), dxdy\big).$$
Moreover, hypothesis \eqref{atzero} assumes the form
$$
\displaystyle\limsup_{t\rightarrow 0}\frac{f(x,t)}{t}< m_0\lambda_1,
$$
where $\lambda_1$ is the first eigenvalue of the problem
$$\left\{\begin{array}{ll}
-\mathcal L_Ku=\lambda u & \mbox{in}\,\,\,\Omega\\
u=0 & \mbox{in}\,\,\,\erre^n\setminus\Omega.
\end{array}\right.$$
Note that a model for $K$ is given by the singular kernel $K(x):=|x|^{-(n+2s)}$ which gives rise to the fractional Laplace operator.
}
\end{remark}

\begin{remark}
\label{finale}\rm{
The previous remark ensures that our result represents an improvement of
\cite[Theorems 1 and 2]{svmountain}, provided
that $M\equiv 1$.}
\end{remark}

{\bf Acknowledgements.}  This paper was written when the first author was visiting professor at the University of Ljubljana in 2014. He expresses his gratitude to the host institution for warm hospitality.
  The manuscript was realized within the auspices of the GNAMPA Project 2014 titled {\it Propriet\`{a} geometriche ed analitiche per problemi non--locali} and the SRA grants P1-0292-0101 and J1-5435-0101.

\medskip


\begin{thebibliography}{99}

\bibitem{8} {\sc C.O. Alves and F.J.S.A. Corr\^{e}a}, \emph{On existence of solutions for a class of problem involving a nonlinear operator},
Comm. Appl. Nonlinear Anal. \textbf{8} (2001), 43-56.

\bibitem{ACM} {\sc C.O. Alves, F.J.S.A. Corr\^{e}a and T.F. Ma},  Positive solutions for a quasilinear
elliptic equation of Kirchhoff type, {\it Comput. Math. Appl.} \textbf{49} (2005), 85-93.

\bibitem{ar} {\sc A. Ambrosetti and P. Rabinowitz},
{\em Dual variational methods in critical point theory and
applications}, {J. Funct. Anal.} \textbf{14} (1973), 349-381.

\bibitem{9} {\sc D. Andrade and T.F. Ma}, \emph{An operator equation suggested by a class of stationary problems}, Comm. Appl.
Nonlinear Anal. \textbf{4} (1997), 65-71.

\bibitem{AP1} {\sc G. Autuori and P. Pucci},  {\em Kirchhoff systems with nonlinear source and boundary
damping terms}, {Commun. Pure Appl. Anal.} \textbf{9} (2010), 1161-1188.

\bibitem{AP2} {\sc G. Autuori and P. Pucci}, {\em Kirchhoff systems with dynamic boundary conditions}, {
Nonlinear Anal.} \textbf{73} (2010), 1952-1965.

\bibitem{AP3} {\sc G. Autuori and P. Pucci}, {\em Local asymptotic stability for polyharmonic Kirchhoff
systems}, {Appl. Anal.} \textbf{90} (2011), 493-514.

\bibitem{AP0} {\sc G. Autuori and P. Pucci}, \emph{Existence of entire solutions for a class of quasilinear elliptic equations}, NoDEA Nonlinear Differential Equations Appl. \textbf{20} (2013), 977-1009.

\bibitem{AP} {\sc G. Autuori and P. Pucci}, \emph{Elliptic problems involving the fractional Laplacian in $\erre^N$}, J. Differential Equations \textbf{255} (2013), 2340-2362.

\bibitem{p3} {\sc G. Autuori, F. Colasuonno and P. Pucci}, {\em On the existence of stationary solutions for higher order $p$-Kirchhoff problems via variational methods}, Commun. Contemp. Math. \textbf{16} (2014), 1450002, pages 43.

\bibitem{brezis} {\sc H. Br\'ezis}, Analyse fonctionelle. Th\'{e}orie et
applications, {\em Masson}, Paris (1983).

\bibitem{caff} {\sc L. Caffarelli, J.M. Roquejoffre and O. Savin}, \emph{Nonlocal minimal surfaces}, Comm. Pure Appl. Math. \textbf{63} (2010), 1111-1144.

\bibitem{10} {\sc M. Chipot and B. Lovat}, \emph{Some remarks on non local elliptic and parabolic problems}, Nonlinear Anal. \textbf{30} (1997), 4619-4627.

\bibitem{11} {\sc M. Chipot and J.F. Rodrigues}, \emph{On a class of nonlocal nonlinear elliptic problems}, RAIRO Mod\'{e}lisation Math.
Anal. Num\'{e}r. \textbf{26} (1992), 447-467.

\bibitem{valpal} {\sc E. Di Nezza, G. Palatucci and E. Valdinoci},
{\em Hitchhiker's guide to the fractional Sobolev spaces},
Bull. Sci. Math. {\bf 136} (2012), no.~5, 521-573.

\bibitem{mio1}	{\sc M. Ferrara, G. Molica Bisci and B. Zhang}, \emph{Existence of weak solutions for non-local fractional problems via Morse theory}, Discrete and Continuous Dynamical Dystems Series B \textbf{19} (2014), 2493-2499.

\bibitem{fv} {\sc A. Fiscella and E. Valdinoci}, \textit{A critical Kirchhoff type problem involving a non-local operator}, Nonlinear Anal. \textbf{94} (2014), 156-170.

\bibitem{k2} {\sc A. Krist\'{a}ly, V. R\u{a}dulescu and Cs. Varga}, {\it Variational Principles in Mathematical Physics, Geometry, and Economics:
Qualitative Analysis of Nonlinear Equations and Unilateral Problems}, Encyclopedia of Mathematics and its Applications \textbf{136}, Cambridge University Press, Cambridge, 2010.

\bibitem{Molica} {\sc G. Molica Bisci}, {\em Fractional equations with bounded primitive},
Appl. Math. Lett. {\textbf{27}} (2014), 53-58.

\bibitem{Molica1} {\sc G. Molica Bisci}, {\em Sequences of weak solutions for fractional equations},
{Math. Res. Lett.} \textbf{21} (2014), 1-13.

\bibitem{Molica2} {\sc G. Molica Bisci and B.A. Pansera}, {\em Three weak solutions for nonlocal fractional equations},
{Adv. Nonlinear Stud.} {\bf 14} (2014), 591-601.

\bibitem{mio2} {\sc G. Molica Bisci and D. Repov\v{s}}, {\em Fractional nonlocal problems involving nonlinearities with bounded primitive}, J. Math. Anal. Appl.  \textbf{420} (2014), 167-176.

\bibitem{Molica3}{\sc G. Molica Bisci and R. Servadei}, {\em A bifurcation result for nonlocal fractional equations}, {Anal. Appl.} {\bf 13} (4) (2015) 371.

    \bibitem{p1} {\sc P. Pucci and Q. Zhang}, {\em Existence of entire solutions for a class of variable exponent elliptic equations}, J. Differential Equations \textbf{257} (2014), 1529-1566.

\bibitem{puccirad} {\sc P. Pucci and V. R\u{a}dulescu}, {\em The impact of the mountain pass theory in
nonlinear analysis: a mathematical survey}, {Boll. Unione Mat. Ital.} Series IX, No. 3, (2010) 543-584.

\bibitem{p2}{\sc P. Pucci and S. Saldi}, {\em Critical stationary Kirchhoff equations in $\erre^n$ involving nonlocal operators}, Rev. Mat. Iberoam. {\bf 32} (1) (2016), 1-22.

\bibitem{rabinowitz} {\sc P.H. Rabinowitz}, Minimax methods in critical point theory with applications to differential equations,  {\em CBMS Reg. Conf. Ser. Math.}, 65, {\em American Mathematical Society}, Providence, RI (1986).

\bibitem{sY}{\sc R. Servadei}, {\em The Yamabe equation in a non-local setting},
Adv. Nonlinear Anal. \textbf{2} (3) (2013), 235-270.

\bibitem{sv}{\sc R. Servadei and E. Valdinoci}, {\em Lewy-Stampacchia type estimates for variational
inequalities driven by nonlocal operators}, Rev. Mat. Iberoam. \textbf{29} (3) (2013), 1091-1126.

\bibitem{svmountain}{\sc R. Servadei and E. Valdinoci}, {\em Mountain Pass solutions
for non-local elliptic operators}, J. Math. Anal. Appl. \textbf{389} (2012), 887-898.

\bibitem{svlinking}{\sc R. Servadei and E. Valdinoci}, {\em Variational methods for non-local operators of elliptic type}, Discrete Contin. Dyn. Syst. \textbf{33}, 5 (2013), 2105-2137.

\bibitem{servadeivaldinociBN}{\sc R. Servadei and E. Valdinoci}, {\em The Br\'{e}zis-Nirenberg result for the fractional Laplacian}, Trans. Amer. Math. Soc. {\bf 367} (2015), 67-102.

\bibitem{struwe} {\sc M. Struwe}, Variational methods,
Applications to nonlinear partial differential equations and Hamiltonian systems,
{\em Ergebnisse der Mathematik und ihrer Grenzgebiete},  3,
{\em Springer Verlag}, Berlin--Heidelberg (1990).

\bibitem{12} {\sc C.F. Vasconcellos},\emph{ On a nonlinear stationary problem in unbounded domains}, Rev. Mat. Univ. Complut.
Madrid \textbf{5} (1992), 309-318.
\end{thebibliography}
\end{document}